\documentclass[12pt,fleqn]{article}
\usepackage{a4}
\usepackage{latexsym,amsmath,amssymb,amsthm,setspace}
\usepackage[english]{babel}
\usepackage{color}

\newtheorem{theorem}{Theorem}
\newtheorem{lemma}{Lemma}

\newtheorem*{conjecture*}{Conjecture}
\newtheorem{corollary}{Corollary}

\newcommand{\Ra}{\Rightarrow}

\newcommand{\llr}{\Longleftrightarrow}
\newcommand{\bs}{\backslash}

\newcommand{\bQ}{\bar{Q}}
\newcommand{\be}{\beta}

\newcommand{\void}{\varnothing}
\newcommand{\eqas}{{\stackrel{\text{a.s.}}{~=~}}} 
\newcommand{\subas}{{\stackrel{\text{a.s.}}{~\subset~}}}

\date{May 16, 2018}

\begin{document}
\setlength{\voffset}{-1.3cm}

\begin{titlepage}
\title{Asymptotic period of an aperiodic Markov chain}

\author{\\ \\Erik A. van Doorn\\
\\Department of Applied Mathematics \\ University of Twente \\
P.O. Box 217, 7500 AE Enschede, The Netherlands\\ E-mail: 
e.a.vandoorn@utwente.nl}

\maketitle
\thispagestyle{empty}

\noindent{\bf Abstract.} 
We introduce the concept of {\em asymptotic period\/} for an irreducible and
aperiodic, discrete-time Markov chain $\mathcal{X}$ on a countable state
space, and develop the theory leading to its formal definition. The asymptotic
period of $\mathcal{X}$ equals one -- its period --  if $\mathcal{X}$ is
recurrent, but may be larger than one if $\mathcal{X}$ is transient;
$\mathcal{X}$ is {\em asymptotically aperiodic\/} if its asymptotic period equals
one. Some sufficient conditions for asymptotic aperiodicity are presented. The
asymptotic period of a birth-death process on the nonnegative integers is
studied in detail and shown to be equal to 1, 2 or $\infty$. Criteria for the
occurrence of each value in terms of the 1-step transition probabilities are
established.



\end{titlepage}

%_______________________________________________________________________________
\section{Introduction}
\label{I}

Let $P:= (P(i,j),~i,j \in S)$ be the matrix of 1-step transition probabilities
of a homogeneous, discrete-time Markov chain $\mathcal{X} :=
\{X(n),~n=0,1,\ldots\}$ on a countably infinite state space $S$, so that the
matrix $P^{(n)} := (P^{(n)}(i,j),~i,j \in S)$ of $n$-step transition
probabilities
\[
P^{(n)}(i,j) := \Pr\{X(m+n)=j\,|\,X(m)=i\}, \quad i,j \in S,
~m,n=0,1,\ldots,
\]
is given by
\[
P^{(n)} = P^n, \quad n=0,1,\ldots.
\]
We will assume throughout that $\mathcal{X}$ is stochastic, irreducible, and
aperiodic. 

Although the Markov chain $\mathcal{X}$ is aperiodic it may happen, if
$\mathcal{X}$ is {\em transient\/}, that in the long run the process evolves
cyclically through a finite number of sets constituting a partition of $S$. This
phenomenon occurs for instance when $\mathcal{X}$ is a transient birth-death
process on the nonnegative integers with only a finite number of positive
self-transition probabilities, for in this case the process will eventually move cyclically
between the even-numbered and the odd-numbered states. It seems natural
then to say that the {\em asymptotic period\/} of $\mathcal{X}$ equals two
or, perhaps, a multiple of two. In our general setting the {\em asymptotic
period\/} of $\mathcal{X}$ may be defined as the maximum number of sets
involved in the type of cyclic behaviour described above. In this paper these ideas
will be formalized, and some of their consequences will be investigated.

After discussing preliminary concepts and results in Section \ref{Pr} we formally
define, in Section \ref{AP}, the {\em asymptotic period\/} of a Markov chain that
is, in a sense to be defined, {\em simple}. Some sufficient conditions for {\em
asymptotic aperiodicity\/} will subsequently be derived. The framework
developed in Section \ref{Pr} draws heavily on the work of Blackwell \cite{B55}
on transient Markov chains, while our definition of asymptotic period resembles
in some aspects the definition of {\em period of an irreducible positive operator\/}
by Moy \cite{M67}, and is directly related to the definition of {\em asymptotic
period of a tail sequence of subsets of\/} $S$, proposed by Abrahamse \cite{A69}
in a setting that is more general than ours. Actually, Abrahamse introduces the
concept of asymptotic period while generalizing Blackwell's results. Our further
elaboration of the concept in a more restricted setting makes it more convenient
for us to build directly on the foundations laid down by Blackwell.

In Section \ref{BD} we investigate asymptotic periodicity in the specific setting
of a birth-death process on the nonnegative integers. We show that the
asymptotic period equals $1$, $2$ or $\infty$, and identify the circumstances
under which each value occurs in terms of the 1-step transition probabilities of
the process. In particular, we establish a necessary and sufficient condition for
asymptotic aperiodicity. 

Our motivation for introducing the concept of {\em asymptotic aperiodicity\/}
has been our aim to gain more insight into the {\em strong ratio limit
property\/}, which is said to prevail if there exist positive constants $R$,
$\mu(i),~i \in S$, and $f(i),~i \in S$, such that
\begin{equation}
\label{srlp}
\lim_{n \to \infty} \frac{P^{(n+m)}(i,j)}{P^{(n)}(k,l)} = R^{-m} \frac
{f(i)\mu(j)}{f(k)\mu(l)}, \quad i,j,k,l \in S,~m \in \mathbb{Z}.
\end{equation}
The strong ratio limit property was enunciated in the setting of recurrent
Markov chains by Orey \cite{O61}, and introduced in the more general setting
at hand by Pruitt \cite{P65}. More recently, Kesten \cite{K95} and Handelmann
\cite{H99} have made substantial contributions, but a satisfactory solution
to the problem of finding conditions for the strong ratio limit property is still
lacking. Since {\em aperiodicity\/} is necessary and sufficient for a
{\em positive recurrent\/} Markov chain to possess the strong ratio limit
property, it is to be expected that {\em asymptotic aperiodicity\/} is a relevant
property in the more general setting at hand. Actually, the relation between
asymptotic aperiodicity and the strong ratio limit property is -- at least in the
general setting -- not clear-cut. However, in a more restricted setting asymptotic
aperiodicity has been shown in \cite{D18} to be sufficient for the strong ratio
limit property. 

We end this introduction with some notation and terminology. Namely, when
$\mathcal{X}$ is a discrete-time birth-death process on the nonnegative
integers -- a process often encountered in what follows -- we write 
\begin{equation}
\label{bdrates}
p_i := P(i,i+1),~~q_{i+1} := P(i+1,i)~~\mbox{and}~~r_i := P(i,i),
\quad i = 0,1,\ldots,
\end{equation}
for the birth, death and self-transition probabilities, respectively. It will
be convenient to define $q_0 := 0$. Since $\mathcal{X}$ is stochastic,
irreducible and aperiodic, we have $p_i>0$,  $q_{i+1}>0$, and $r_i\geq 0$ for
$i \geq 0$, with $r_i>0$ for at least one state $i$, while $p_i+q_i+r_i=1$ for
$i\geq 0$. In what follows a {\em birth-death process\/} will always refer to
a discrete-time birth-death process on the nonnegative integers.

%______________________________________________________________________________
\section{Preliminaries}
\label{Pr}

We start off by introducing some further notation and terminology related to
the Markov chain $\mathcal{X} = \{X(n),~n=0,1,\ldots\}$. 
By $\mathbb{P}$ we denote the probability measure on the set of sample paths
induced by $P$ and the (unspecified) initial distribution. Recall that a nonzero
function $f$ on $S$ is called a {\em harmonic function\/} (or {\em invariant
vector\/}) for $P$ (or, for $\mathcal{X}$) if
\begin{equation}
\label{f}
Pf(i) := \sum_{j \in S} P(i,j)f(j) = f(i), \quad i \in S.
\end{equation}
Evidently, in our setting the constant function is a harmonic function for $P$.

For $C \subset S$ we define the events
\[
U(C) := \cap_{n=0}^\infty\cup_{k=n}^\infty\{X(k)\in C\}
\mbox{~~and~~}
L(C) := \cup_{n=0}^\infty\cap_{k=n}^\infty\{X(k)\in C\},
\]
and we let
\[
\mathcal{T} := \{C \subset S\,|\, U(C) \eqas \void\},
\]
that is, $C \in \mathcal{T}$ if $\mathbb{P}(X(n) \in C$ infinitely often$) = 0$,
and
\[
\mathcal{R} := \{C \subset S\,|\, U(C) \eqas L(C)\},
\]
that is, $C \in \mathcal{R}$ if the events \{$X(n) \in C$ infinitely often\} and
\{$X(n) \in C$ for $n$ sufficiently large\} are almost surely equal. In the
terminology of Revuz \cite[Sect.~2.3]{R75} $\mathcal{T}$ is the collection
of {\em transient\/} sets and $\mathcal{R}$ is the collection of {\em regular\/}
sets. Evidently, $\mathcal{T} \subset \mathcal{R}$, while it is not difficult
to see that $\mathcal{R}$ is closed under finite union and complementation,
and hence a field.  Note that $\mathcal{T}$ and $\mathcal{R}$ are independent
of the initial distribution, since, by the irreducibility of $\mathcal{X}$,
$\mathbb{P}(U(C))$ and $\mathbb{P}(U(C)\bs L(C))$ are zero or positive for all
initial states (and hence all initial distributions) simultaneously.

We will say that two regular sets $C_1$ and $C_2$ are {\em equivalent\/} if
their symmetric difference $C_1 \Delta C_2 := (C_1 \cup C_2)\bs(C_1 \cap C_2)$
is transient, and {\em almost disjoint\/} if their intersection $C_1 \cap C_2$
is transient. Following Blackwell \cite{B55} (see also Chung \cite[Section
I.17]{C67}), we call a subset $C \subset S$ {\em almost closed\/} if $C \notin
\mathcal{T}$ and $C \in \mathcal{R}$. An almost closed set $C$ is said to be
{\em atomic\/} if $C$ does not contain two disjoint almost closed subsets.
The relevance of these concepts comes to light in the next theorem.

\begin{theorem}
\label{blackwell1}
{\em (Blackwell \cite{B55})}
Associated with the Markov chain $\mathcal{X}$ is a finite or countable
collection $\{C_1,C_2,\ldots\}$ of disjoint almost closed sets, which is
unique up to equivalence and such that\\
(i) every $C_i$, except at most one, is atomic;\\
(ii) the nonatomic $C_i$, if present, contains no atomic subsets and
consists of transient states;\\
(iii) $\sum_i \mathbb{P}(L(C_i)) = 1$.
\end{theorem}

\noindent
A collection of sets $\{C_1,C_2,\ldots\}$ satisfying the conditions in the
theorem will be called a {\em Blackwell decomposition\/} (of $S$) for
$\mathcal{X}$. A set $C \subset S$ is a {\em Blackwell component\/} (of $S$)
for $\mathcal{X}$ if there exists a Blackwell decomposition for $\mathcal{X}$
such that $C$ is one of the almost closed sets in the decomposition. The
uniqueness up to equivalence of the Blackwell decomposition for $\mathcal{X}$
implies that if $C_1$ and $C_2$ are Blackwell components, then they are either
equivalent or almost disjoint. The number of almost closed sets in the
Blackwell decomposition for $\mathcal{X}$ will be denoted by
$\beta(\mathcal{X})$.
If $\beta(\mathcal{X}) = 1$ then $\mathcal{X}$ is called {\em simple\/}, and
a simple Markov chain is called {\em atomic\/} or {\em nonatomic\/} according to
the type of its state space. Evidently, if $\mathcal{X}$ is simple and
nonatomic then $S$ does not contain atomic subsets, but infinitely many disjoint
almost closed subsets.
%For $S$, being nonatomic, contains two disjoint almost closed subsets, each
%of which is nonatomic and hence contains two disjoint almost closed subsets, etc.
It will be useful to observe the following.

\begin{lemma}
\label{bounded}
Let $S = \{0,1,\dots\}$ and $\mathcal{X}$ have jumps that are uniformly
bounded by $M$. Then $\beta(\mathcal{X}) \leq M$, and every Blackwell component
for $\mathcal{X}$ is atomic.
\end{lemma}
\noindent
{\bf Proof.}~
Let $C$ be an almost closed set for $\mathcal{X}$ and let $s_1<s_2<\dots$
denote the states of $C$. We claim that there exists a constant $N$ such that
for every $n\geq N$ we have $s_{n+1} \leq s_n+M$. Indeed, if $s_{n+1}>s_n+M$,
then the process will leave $C$ when it leaves the set $\{s_1,s_2,\dots,s_n\}$.
The irreducibility of $S$ insures that a visit to this finite set of states will
almost surely be followed by a departure from the set. So if, for each $N$,
there is an integer $n\geq N$ such that $s_{n+1}>s_n+M$, then each entrance
in $C$ is almost surely followed by a departure from $C$, and hence
$\mathbb{P}(L(C))=0$, contradicting the fact that $C$ is almost closed.

Next, let $\{C_1,C_2,\dots,C_\beta\}$, with $\beta\equiv\beta(\mathcal{X})$,
be a Blackwell decomposition for $\mathcal{X}$, $s_1^{(i)}<s_2^{(i)}<\dots$ 
the states of $C_i$, and $N_i$ such that for every $n\geq N_i$ we have
$s_{n+1}^{(i)} \leq s_n^{(i)}+M$. If $\beta > M$, then, choosing
\[
s = \max_{1\leq i\leq M+1}\{s_{N_i}^{(i)}\},
\]
the set $\{s+1,s+2,\dots,s+M\}$ must have a nonempty intersection with each
of the disjoint sets $C_1,C_2,\dots,C_{M+1}$, which is clearly impossible. Hence,
$\beta\leq M$.

Finally, let $C$ be a Blackwell component for $\mathcal{X}$ and suppose $C$ is
nonatomic. Then $C$ contains infinitely many disjoint almost closed subsets,
so we can choose $M+1$ disjoint almost closed subsets $C_1,C_2,\dots,C_{M+1}$
of $C$. By the same argument as before there must be a state $s$ in $C$ such
that each of the disjoint sets $C_1,C_2,\dots,C_{M+1}$ shares a state with the
set $\{s+1,s+2,\dots,s+M\}$. This is impossible, so $C$ must be atomic. $\Box$
\bigskip

\noindent
A criterion for deciding whether a Markov chain is simple and atomic is given in
the next theorem. 

\begin{theorem}
\label{blackwell2}
{\em (Blackwell \cite{B55})}
The Markov chain $\mathcal{X}$ is simple and atomic if and only if the only
bounded harmonic function for $\mathcal{X}$ is the constant function.
\end{theorem}

\noindent
As an aside we note that when $\mathcal{X}$ is transient -- the setting of
primary interest to us -- and the constant function is the only bounded
harmonic function, then there is precisely one escape route to infinity, or, in
the terminology of Hou and Guo \cite{HG88} (see, in particular, Sections 7.13
and 7.16), the {\em exit space\/} of $\mathcal{X}$ contains exactly one
{\em atomic exit point\/}.

Of course, the existence, up to a multiplicative
constant, of a unique {\em bounded\/} harmonic function does not, in
general, preclude the existence of an {\em unbounded\/} harmonic function.
%This happens, for instance, in the setting of a birth-death process on the
%integers when the process has a drift in precisely one direction.
But when $\mathcal{X}$ is recurrent the constant function happens to be
the only (bounded or unbounded) harmonic function (see, for example,
Chung \cite[Theorem I.7.6]{C67}). It follows in particular that $\mathcal{X}$
is simple and atomic if $\mathcal{X}$ is recurrent.

A function $f$ on the space $\Omega := \{(\omega_0,\omega_1,\ldots)\,|\,
\omega_i \in S,~i=0,1,\ldots\}$ will be called $m$-{\em invariant\/} if, for
every $\omega := (\omega_0,\omega_1,\ldots) \in
\Omega,~f(\omega)=f(\theta^m\omega),$ where $\theta$ is the shift operator
$\theta(\omega_0, \omega_1,\ldots)=(\omega_1,\omega_2,\ldots),$ and 
$\theta^m\omega = \theta(\theta^{m-1}\omega).$ We also use
the notation $\theta^m E := \{\theta^m\omega\,|\,\omega \in E\}$, for
$E \subset \Omega$. An event is called 
$m$-{\em invariant\/} if its indicator function is $m$-invariant. A
$1$-invariant event is simply referred to as {\em invariant\/}.
Evidently, the collection of invariant events constitutes a $\sigma$-field.
We shall need another result of Blackwell's, involving invariant events
(see \cite[Theorem 5]{A69} for a generalization).

\begin{theorem}
\label{blackwell3}
{\em (Blackwell \cite{B55})}
For any invariant event $E$ there is a $C \in \mathcal{R}$ such that
$E \eqas U(C)$.
\end{theorem}

\noindent
Note that the event $U(C)$ is actually invariant for {\em any\/} subset $C$
of $S$, so for every $C \subset S$ there must be a {\em regular\/} set
$\tilde{C}$ such that $U(C) \eqas U(\tilde{C})$. It follows in particular that
every invariant event has probability zero or one if $\mathcal{X}$ is simple
and atomic.

The regular set corresponding to an invariant event is unique up to
equivalence. For if $C_1$ and $C_2$ are regular sets satisfying
$U(C_1) \eqas U(C_2)$, then
\[
U(C_1\bs C_2) \subset U(C_1)\bs L(C_2) \eqas U(C_2)\bs L(C_2) \eqas \void,
\]
and similarly with $C_1$ and $C_2$ interchanged. Since
$U(C_1\Delta C_2) \subset U(C_1\bs C_2) \cup U(C_2\bs C_1)$,
it follows that $C_1 \Delta C_2$ must be transient.
So, up to events of probability zero, the $\sigma$-field of invariant events
is identical with the $\sigma$-field of events of the form $U(C)$ with
$C \in \mathcal{R}$. 

Theorem \ref{blackwell3} plays a crucial role in the proof of Theorem
\ref{prop}, which involves
$\mathcal{X}^{(m)} := \{X^{(m)}(n) \equiv X(nm),~n=0,1,\ldots\}$,
the $m$-step Markov chain associated with $\mathcal{X}$, and is
instrumental in our definition of asymptotic period. For $C \subset S$ we let
\[
U^{(m)}(C) := \cap_{n=0}^\infty\cup_{k=n}^\infty\{X(km)\in C\}
\mbox{~and~}
L^{(m)}(C) := \cup_{n=0}^\infty\cap_{k=n}^\infty\{X(km)\in C\},
\]
so that $U^{(1)}(C) = U(C)$ and $L^{(1)}(C) = L(C)$. Since $E = \theta^m E$
if (and only if) $E$ is $m$-invariant, we have $\theta^mU^{(m)}(C) = U^{(m)}(C)$.
But actually we have, more generally,
\begin{equation}
\label{UmC}
\theta^{m-j}U^{(m)}(C) = \cap_{n=0}^\infty\cup_{k=n}^\infty\{X(km+j)\in C\},
\quad j=0,1,\dots,m-1,
\end{equation}
and
\begin{equation}
\label{LmC}
\theta^{m-j}L^{(m)}(C) = \cup_{n=0}^\infty\cap_{k=n}^\infty\{X(km+j)\in C\},
\quad j=0,1,\dots,m-1,
\end{equation}
as can easily be verified. The following simple observation will prove useful.
\begin{lemma}
\label{aux1}
Let $E$ be an $m$-invariant event for some $m\geq 1$. Then, for all $i \geq 1$,
\[
\quad E \eqas \void ~\llr~ \theta^i E \eqas \void.
\]
\end{lemma}
\noindent
{\bf Proof.}~
If $\mathbb{P}(E)>0$ there must be a state $s$,
say, such that $\mathbb{P}(E\,|\,X(0)=s)>0$. Moreover, aperiodicity
and irreducibility of the chain imply that there is an integer $k$ such that
$\mathbb{P}(X(k m-i)=s) > 0$. Since, by Theorem \ref{blackwell3},
$E \eqas U^{(m)}(C)$ for some set $C$, we obviously have 
$\mathbb{P}(\theta^i E\,|\,X(km-i)=s) = \mathbb{P}(E\,|\,X(0)=s)$.
Hence, if $\mathbb{P}(E)>0$, then
\[
\begin{array}{l@{}l}
\mathbb{P}(\theta^i E)~ 
&\geq \mathbb{P}(\theta^i E\,|\,X(k m-i)=s)\mathbb{P}(X(k m-i)=s)
\vspace{0.2cm}\\
&= \mathbb{P}(E\,|\,X(0)=s)\mathbb{P}(X(k m-i)=s) > 0.
\end{array}
\]
The same argument with $E$ and $\theta^i E$ interchanged and $km-i$
replaced by $km+i$ yields the converse. $\Box$
\bigskip

\noindent
Before stating and proving Theorem \ref{prop} we establish some
additional auxiliary lemmas. In what follows we write $E\subas F$ for
$E\bs F \eqas \void$.

\begin{lemma}
\label{aux2}
Let $E_1$ and $E_2$ be $m$-invariant events  for some $m\geq 1$.
Then, for all $i \geq 0$ and $j\geq 0$, \vspace{0.2cm}\\
(i) $\quad  E_1 \subas \theta^j E_2  ~\llr~ \theta^i E_1 \subas \theta^{i+j} E_2$,
\vspace{0.2cm}\\
(ii) $\quad E_1 \eqas \theta^j E_2  ~\llr~ \theta^i E_1 \eqas \theta^{i+j} E_2$.
\end{lemma}
\noindent
{\bf Proof.}~
The event $E_1\bs \theta^j E_2$ is $m$-invariant, so, by Lemma \ref{aux1},
we have
\[
E_1\bs \theta^j E_2 \eqas \void ~~\llr~~ \theta^i (E_1\bs \theta^j E_2) \eqas \void,
\]
which implies the first statement. Moreover, the first statement remains
valid, by a similar argument, if we interchange the sets $E_1$ and
$\theta^j E_2$. Combining both results yields the second statement. $\Box$
\bigskip

\noindent
Note that the second statement of this lemma generalizes Lemma \ref{aux1}.
The next auxiliary result is a straightforward corollary of the previous lemma.

\begin{lemma}
\label{aux3}
Let $E$ be an $m$-invariant event  for some $m\geq 1$.
Then, for all $j \geq 0$ and $k_2 \geq k_1 \geq 0$,\vspace{0.2cm}\\
(i) $\quad E \subas \theta^j E ~\Ra~
\theta^{k_1j}E \subas \theta^{k_2j} E$,
\vspace{0.2cm}\\
(ii) $\quad E \eqas \theta^j E ~\Ra~
\theta^{k_1j}E \eqas \theta^{k_2j} E$.
\end{lemma}

\noindent
Our final preparatory lemma is the following.

\begin{lemma}
\label{aux4}
Let $C_1$ and $C_2$ be subsets of $S$ that are regular with respect to
$\mathcal{X}^{(m)}$ for some $m\geq 1$. Then
\[
U^{(m)}(C_1 \cap C_2) \eqas U^{(m)}(C_1) \cap U^{(m)}(C_2).
\]
\end{lemma}
\noindent
{\bf Proof.}~
We clearly have
\[
U^{(m)}(C_1 \cap C_2) \subset U^{(m)}(C_1) \cap U^{(m)}(C_2) \eqas
L^{(m)}(C_1) \cap L^{(m)}(C_2).
\]
Since
\[
L^{(m)}(C_1) \cap L^{(m)}(C_2) = L^{(m)}(C_1 \cap C_2) \subset
U^{(m)}(C_1 \cap C_2),
\]
the result follows. $\Box$

\begin{theorem}
\label{prop}
If $\mathcal{X}$ is simple and atomic and $m>1$, then $\beta \equiv
\beta(\mathcal{X}^{(m)})$ is a divisor of $m$ and the Blackwell decomposition
for $\mathcal{X}^{(m)}$ consists of a collection
$\{C_0,C_1,\ldots,C_{\beta-1}\}$ of disjoint atomic almost closed sets, which
can be chosen such that, for each $i = 0,1,\ldots,\beta-1$,
\begin{equation}
\label{structure0}
\mathbb{P}(\theta^j L^{(m)}(C_{i+1\,(\mbox{\scriptsize{mod}}\,\beta)})\,|\,
\theta^{j+1} L^{(m)}(C_i)) = 1, \quad j = 0,1,\ldots,m-1.
\end{equation}
If $\mathcal{X}$ is simple and nonatomic, then $\mathcal{X}^{(m)}$ is simple
and nonatomic for all $m\geq 1$.
\end{theorem}
%The assumption that $\mathcal{X}$ be simple is essential here, since the
%Blackwell components for $\mathcal{X}$, if there are at least two, can have
%very different structures. 

\noindent
{\bf Proof.}~
First suppose $\mathcal{X}$ is simple and atomic. Let $C_0$ be a Blackwell
component for $\mathcal{X}^{(m)}$ and assume, for the time being, that $C_0$
is atomic.
Since $\theta^iU^{(m)}(C_0)$ is $m$-invariant for all $i$, we can apply Theorem
\ref{blackwell3} to $\mathcal{X}^{(m)}$ and conclude that there is a sequence
$C_1,C_2,\ldots$ of regular sets (with respect to $\mathcal{X}^{(m)}$) such that 
\begin{equation}
\label{equiv}
\theta^i U^{(m)}(C_0) \eqas U^{(m)}(C_i), \quad i=1,2,\ldots.
\end{equation}
By Lemma \ref{aux1} the sets $C_i$ are almost closed, since $C_0$ is
almost closed. Also, by Lemma \ref{aux2}, we have
\[
U^{(m)}(C_{i+1})\eqas \theta^{i+1} U^{(m)}(C_0) \eqas \theta U^{(m)}(C_i),
\]
and hence
\begin{equation}
\label{structure1}
L^{(m)}(C_{i+1})\eqas \theta L^{(m)}(C_i), \quad i=0,1,2,\dots.
\end{equation}

Next defining
\begin{equation}
\label{b}
b := \min\{i\geq 1\,|\,\theta^i U^{(m)}(C_0) \eqas U^{(m)}(C_0)\},
\end{equation}
we have $b \leq m$ since $U^{(m)}(C_0)$ is $m$-invariant. Also, $b$ must
be a divisor of $m$, for otherwise, by Lemma \ref{aux3}, we would have
\[
U^{(m)}(C_0) \eqas \theta^{\ell b} U^{(m)}(C_0) =
\theta^{m+i} U^{(m)}(C_0) =  \theta^i U^{(m)}(C_0),
\]
with $\ell=\min\{k\in \mathbb{N}\,|\,kb>m\}$ and $i=\ell b-m<b$, contradicting
\eqref{b}.
For $i \geq b$ we have, by Lemma \ref{aux2},
\[
U^{(m)}(C_i) \eqas \theta^i U^{(m)}(C_0) \eqas
\theta^{i-b} U^{(m)}(C_0) \eqas U^{(m)}(C_{i-b}),
\]
so that $C_i$ and $C_{i-b}$ are equivalent (with respect to $\mathcal{X}^{(m)}$).
We can therefore replace \eqref{structure1} by
\[
L^{(m)}(C_{i+1\,(\mbox{\scriptsize{mod}}\,b)})\eqas
\theta L^{(m)}(C_i), \quad i=0,1,2,\dots,b-1.
\]
But in view of Lemma \ref{aux2} we can actually write, for any value of $j$,
\begin{equation}
\label{structure2}
\theta^j L^{(m)}(C_{i+1\,(\mbox{\scriptsize{mod}}\,b)})\eqas
\theta^{j+1} L^{(m)}(C_i), \quad i=0,1,2,\dots,b-1,
\end{equation}
so that \eqref{structure0} prevails for $i = 0,1,\ldots,b-1$.

Our next step will be to prove that the sets $C_0,C_1,\ldots,C_{b-1}$ are
almost disjoint. Since the collection of sets that are regular with respect to
$\mathcal{X}^{(m)}$ constitutes a field, the sets $C_0\bs C_i$ and
$C_0\cap C_i$, with $0<i<b$, are regular. But $C_0$, being an atomic Blackwell
component for $\mathcal{X}^{(m)}$, cannot contain two almost closed subsets,
so that either $C_0\bs C_i$ or $C_0\cap C_i$ must be transient.
If $C_0\bs C_i$ is transient, then 
\[
U^{(m)}(C_0)\bs U^{(m)}(C_i) \subset U^{(m)}(C_0\bs C_i) \eqas \void,
\]
which implies 
\[
U^{(m)}(C_0) \eqas U^{(m)}(C_0)\cap U^{(m)}(C_i) \subset
U^{(m)}(C_i) \eqas \theta^i U^{(m)}(C_0),
\]
that is, $U^{(m)}(C_0) \subas \theta^i U^{(m)}(C_0)$. But then, by Lemma
\ref{aux3},
\[
\theta^i U^{(m)}(C_0) \subas \theta^{bi} U^{(m)}(C_0)
\eqas U^{(m)}(C_0),
\]
so that $U^{(m)}(C_0) \eqas \theta^i U^{(m)}(C_0)$, contradicting \eqref{b}.
So we conclude, for $0<i<b$, that $C_0\cap C_i$ is transient, and hence that
$C_0$ and $C_i$, are almost disjoint. It subsequently follows that $C_i$ and
$C_j$, with $0 \leq i < j < b$, are also almost disjoint. Indeed, $C_0$ and
$C_{j-i}$ being almost disjoint, we have, by Lemma \ref{aux4},
\[
U^{(m)}(C_0) \cap \theta^{j-i}U^{(m)}(C_0) \eqas
U^{(m)}(C_0) \cap U^{(m)}(C_{j-i}) \eqas U^{(m)}(C_0 \cap C_{j-i}) \eqas \void.
\]
Hence, by Lemma \ref{aux1} and Lemma \ref{aux4},
\[
U^{(m)}(C_i \cap C_j) \eqas U^{(m)}(C_i) \cap U^{(m)}(C_j) \eqas
\theta^i \left(U^{(m)}(C_0) \cap \theta^{j-i}U^{(m)}(C_0) \right) \eqas \void,
\]
establishing our claim. It is no restriction of generality to assume that the sets
$C_0,C_1,\ldots,C_{b-1}$ are actually {\em disjoint\/} (rather than {\em almost
disjoint\/}), since replacing $C_i$ by the equivalent set $C'_i$, where
$C'_0 = C_0$ and $C'_i = C_i\bs\cup_{j<i}C_j,~ i = 1,\dots,b-1$, does not
disturb the validity of \eqref{equiv}.

Our next step will be to show that $\{C_0,C_1,\ldots,C_{b-1}\}$ constitutes a
Blackwell decomposition for $\mathcal{X}^{(m)}$, still assuming the Blackwell
component $C_0$ to be atomic. First note that $\cup_{i=0}^{b-1}C_i$ is regular
with respect to $\mathcal{X}$. Indeed, by definition of $C_i$ and in view of
\eqref{UmC} and \eqref{LmC}, we have
\[
U(\cup_{i=0}^{b-1}C_i) = \cup_{i=0}^{b-1}\cup_{j=0}^{m-1}\theta^{m-j}U^{(m)}(C_i)
\eqas \cup_{i=0}^{b-1}\cup_{j=0}^{m-1}\theta^{m-j}L^{(m)}(C_i)
\subas 
L(\cup_{i=0}^{b-1}C_i),
\]
so that $U(\cup_{i=0}^{b-1}C_i) \eqas L(\cup_{i=0}^{b-1}C_i)$. Moreover, 
\[
\mathbb{P}(U(\cup_{i=0}^{b-1}C_i)) \geq
\mathbb{P}(U^{(m)}(\cup_{i=0}^{b-1}C_i)) \geq \mathbb{P}(U^{(m)}(C_0)) > 0,
\]
so $\cup_{i=0}^{b-1} C_i$ is in fact almost closed. It follows, $\mathcal{X}$
being simple and atomic, that $\cup_{i=0}^{b-1} C_i$ and $S$ are equivalent
with respect to $\mathcal{X}$. As a consequence
\[
\mathbb{P}(U^{(m)}(S\backslash\cup_{i=0}^{b-1}C_i))
\le \mathbb{P}(U(S\backslash\cup_{i=0}^{b-1}C_i))=0,
\]
that is, $\cup_{i=0}^{b-1} C_i$ and $S$ are also equivalent with respect to
$\mathcal{X}^{(m)}$. Hence
\[ 
\sum_{i=0}^{b-1} \mathbb{P}(U^{(m)}(C_i))
\ge\mathbb{P}(U^{(m)}(\cup_{i=0}^{b-1}C_i))
\ge 1-\mathbb{P}(U^{(m)}(S\backslash\cup_{i=0}^{b-1}C_i))=1,
\]
so that
\[
\sum_{i=0}^{b-1} \mathbb{P}(L^{(m)}(C_i)) = 
\sum_{i=0}^{b-1} \mathbb{P}(U^{(m)}(C_i)) = 1.
\]
If $b=1$ then $C_0$ and $S$ are equivalent with respect to $\mathcal{X}^{(m)}$,
so that $\beta(\mathcal{X}^{(m)})=1$, and we are done.
So suppose $b>1$ and let $\Gamma$ be an arbitrary almost closed subset of
$C_i,~0<i<b$.
Since $\theta^{b-i}U^{(m)}(\Gamma)$ is invariant with respect to
$\mathcal{X}^{(m)}$, there exists, by Theorem \ref{blackwell3}, a regular set
$\Gamma_0$ such that $\theta^{b-i}U^{(m)}(\Gamma) \eqas U^{(m)}(\Gamma_0)$.
Lemma \ref{aux1} implies that $\Gamma_0$ is almost closed, while,
by \eqref{structure2}, $U^{(m)}(\Gamma_0) \subas U^{(m)}(C_0)$. But
since $C_0$ is atomic, we must actually have $U^{(m)}(\Gamma_0) \eqas U^{(m)}(C_0)$.
Hence, by Lemma \ref{aux2},
\[
U^{(m)}(\Gamma) = \theta^i(\theta^{b-i} U^{(m)}(\Gamma)) \eqas
\theta^i U^{(m)}(\Gamma_0) \eqas \theta^i U^{(m)}(C_0) \eqas U^{(m)}(C_i),
\]
so that $\Gamma$ and $C_i$ are equivalent. Hence $C_i$ is atomic.
So we conclude that if $C_0$ is atomic then $\{C_0,C_1,\ldots,C_{b-1}\}$
constitutes a Blackwell decomposition for $\mathcal{X}^{(m)}$ (with atomic
components) and hence $\beta\equiv\beta(\mathcal{X}^{(m)}) = b$, a divisor of $m$.

We will now show that, in fact, each component in the Blackwell decomposition
for $\mathcal{X}^{(m)}$ has to be atomic if $\mathcal{X}$ is simple and atomic.
If $\beta(\mathcal{X}^{(m)})>1$, we could replace $C_0$ in the preceding
argument by an atomic Blackwell component for $\mathcal{X}^{(m)}$, and
subsequently reach a contradiction, since all the components in the
Blackwell decomposition for $\mathcal{X}^{(m)}$ have to be atomic if $C_0$
is atomic. So it remains to consider the case $\beta(\mathcal{X}^{(m)})=1$.
Assuming $S$ to be nonatomic with respect to $\mathcal{X}^{(m)}$, there
are almost closed sets that are {\em not\/} equivalent to $S$. Let
$\Gamma_0$ be such a set. Then, by Theorem \ref{blackwell3}, there are
sets $\Gamma_i$, regular with respect to $\mathcal{X}^{(m)}$
and unique up to equivalence, such that
\[
\theta^i U^{(m)}(\Gamma_0) \eqas U^{(m)}(\Gamma_i), \quad i=1,2,\ldots.
\]
Copying the argument following \eqref{equiv} up to and including
\eqref{structure2} with $C_i$ replaced by $\Gamma_i$, we conclude from the
analogue of \eqref{structure2} that $\cup_{i=0}^{b-1}\Gamma_i$ is regular
with respect to $\mathcal{X}$, while 
\[
\mathbb{P}(U(\cup_{i=0}^{b-1}\Gamma_i)) \geq
\mathbb{P}(U^{(m)}(\cup_{i=0}^{b-1}\Gamma_i)) \geq
\mathbb{P}(U^{(m)}(\Gamma_0)) > 0.
\]
So $\cup_{i=0}^{b-1} \Gamma_i$ is in fact almost closed, and it follows,
$\mathcal{X}$ being simple and atomic, that $\cup_{i=0}^{b-1} \Gamma_i$
and $S$ are equivalent with respect to $\mathcal{X}$.

It is no restriction of generality to assume that the sets
$\Gamma_0,\Gamma_1,\ldots,\Gamma_{b-1}$ are disjoint. Indeed,
$\Gamma_0\bs\Gamma_i$ cannot be transient, by the same argument we have
used earlier for $C_0\bs C_i$. Hence, the collection of regular sets constituting
a field, $\Gamma_0\bs\Gamma_i$ must be almost closed with respect to
$\mathcal{X}^{(m)}$. So, if $\Gamma_0\cap\Gamma_i$ is
{\em not\/} transient, we may replace $\Gamma_0$ by $\Gamma_0\bs\Gamma_i$
in the preceding argument and end up with new sets
$\Gamma_0,\Gamma_1,\ldots,\Gamma_{b-1}$ such that 
$\Gamma_0\cap\Gamma_i$ is transient. Repeating the procedure if necessary,
we reach, after less than $b$ steps, a situation in which
$\Gamma_0\cap\Gamma_i$ is transient for each $i<b$. It follows, by the same
argument we have used before for the $C_i$'s, that {\em all\/} $\Gamma_i$'s
are almost disjoint and by a similar adaptation as before for the $C_i$'s we can
actually make them disjoint without essentially changing the situation.
But if $\Gamma_0,\Gamma_1,\ldots,\Gamma_{b-1}$ are disjoint almost closed
sets such that the analogue of \eqref{structure2} is satisfied, and
$\cup_{i=0}^{b-1} \Gamma_i$ and $S$ are equivalent with respect to
$\mathcal{X}$, then $\{\Gamma_0,\Gamma_1,\ldots,\Gamma_{b-1}\}$
constitutes a Blackwell decomposition for $\mathcal{X}^{(m)}$, which, since
$\beta(\mathcal{X}^{(m)})=1$, implies $b=1$, and hence that $\Gamma_0$
and $S$ are equivalent, contradicting our assumption on $\Gamma_0$. So if
$\mathcal{X}$ is simple and atomic and $\beta(\mathcal{X}^{(m)})=1$, then
$S$ has to be atomic. Summarizing we conclude that every component in the
Blackwell decomposition of $S$ for $\mathcal{X}^{(m)}$ must be atomic if
$\mathcal{X}$ is simple and atomic.

Finally, suppose $\mathcal{X}$ is simple and nonatomic. Evidently, each subset
of $S$ that is almost closed with respect to $\mathcal{X}$ contains a subset
that is almost closed with respect to $\mathcal{X}^{(m)}$, and it follows that a
nonatomic almost closed set with respect to $\mathcal{X}$ must contain a
nonatomic almost closed set with respect to $\mathcal{X}^{(m)}$. So $S$ must
contain a nonatomic almost closed set with respect to $\mathcal{X}^{(m)}$.
We have seen that all components in the Blackwell decomposition of $S$ for
$\mathcal{X}^{(m)}$ must be atomic if $\beta(\mathcal{X}^{(m)})>1$, so the
only remaining possibility is that $\mathcal{X}^{(m)}$ is simple and nonatomic. 
$\Box$
\bigskip

\noindent
Note that \eqref{structure0} is equivalent to stating that
for $j=0,1,\dots,m-1$,
\[
\begin{aligned}
&\{X(km+j) \in C_i ~\mbox{for $k$ sufficiently large}\}\\
&\hspace{3cm} \eqas\{X(km+j+1) \in
C_{i+1\,(\mbox{\scriptsize{mod}}\,\beta)}~\mbox{for $k$ sufficiently large}\}\\
\end{aligned}
\]
In what follows we will refer to a Blackwell decomposition of $S$ for
$\mathcal{X}^{(m)}$ with this property as a {\em cyclic\/} decomposition.

Theorem \ref{prop} provides the framework for the formal definition of the
asymptotic period of a simple Markov chain in the next section. We conclude this
section with a series of lemmas and corollaries, which supply further information
on $\beta(\mathcal{X}^{(m)})$.

\begin{lemma}
\label{multiple}
Let $\mathcal{X}$ be simple and atomic, and $m\geq 1$. Then a Blackwell
component for $\mathcal{X}^{(m)}$ is almost closed with respect to
$\mathcal{X}^{(k\beta)}$ for all $k\geq 1$, where
$\beta \equiv \beta(\mathcal{X}^{(m)})$. Also,
$\beta(\mathcal{X}^{(\beta)}) = \beta$.
\end{lemma}

\noindent
{\bf Proof.}~
Let $C$ be a Blackwell component for $\mathcal{X}^{(m)}$. As a consequence of
\eqref{structure0} we have $U^{(\beta)}(C) \eqas L^{(\beta)}(C)$, and hence
$U^{(k\beta)}(C) \eqas L^{(k\beta)}(C)$ for any $k\geq 1$. Also,
\[
\mathbb{P}(L^{(k\beta)}(C)) \geq \mathbb{P}(L^{(\beta)}(C)) =
\mathbb{P}(U^{(\beta)}(C)) \geq  \mathbb{P}(U^{(m)}(C)) > 0,
\]
since $\beta$ is a divisor of $m$. So we conclude that $C$ is almost closed
with respect to $\mathcal{X}^{(k\beta)}$. It follows in particular that a
Blackwell component for $\mathcal{X}^{(m)}$ must contain a Blackwell component
for $\mathcal{X}^{({\beta})}$. Hence $\beta(\mathcal{X}^{(\beta)}) \geq \beta$,
and so $\beta(\mathcal{X}^{(\beta)}) = \beta$, since
$\beta(\mathcal{X}^{(\beta)})$ is a divisor of $\beta$.  $\Box$
\bigskip

\noindent
The following corollary is immediate.
\begin{corollary}
\label{cor1}
Let $\mathcal{X}$ be simple. If $\beta(\mathcal{X}^{(m)}) < m$ for all
$m>1$, then $\beta(\mathcal{X}^{(m)}) = 1$ for all $m$.
\end{corollary}

\begin{lemma}
\label{divisor}
Let $\mathcal{X}$ be simple and $k,\ell\geq 1$. Then
$\beta(\mathcal{X}^{(k\ell)}) = \kappa \beta(\mathcal{X}^{(\ell)})$, where
$\kappa$ is a divisor of $\beta(\mathcal{X}^{(k)})$.
\end{lemma}

\noindent
{\bf Proof.}~
If $\mathcal{X}$ is nonatomic then, by Theorem \ref{prop},
$\beta(\mathcal{X}^{(m)}) = 1$ for all $m$, so that the statement is trivially
true. So let us assume that $\mathcal{X}$ is simple and atomic. We write
$\beta_\ell \equiv \beta(\mathcal{X}^{(\ell)})$, and  denote the (atomic)
Blackwell components for $\mathcal{X}^{(\ell)}$ by
$B_0,B_1,\dots,B_{\beta_{\ell}}$. By the previous lemma these sets are almost
closed with respect to $\mathcal{X}^{(k\ell)}$, so each $B_i$ must contain at
least one Blackwell component for $\mathcal{X}^{(k\ell)}$. Let $C_0 \subset
B_0$ be such a Blackwell component and consider the sets $C_i$ defined in the
proof of Theorem \ref{prop} in terms of $C_0$ and $m=k\ell$. We let
\[
\kappa := \min\{k\geq 1\,|\,\theta^{k\beta_{\ell}} U^{(k\ell)}(C_0)
\eqas U^{(k\ell)}(C_0)\},
\]
and claim that $\kappa\beta_\ell = \beta(\mathcal{X}^{(k\ell)})$. 

To prove the claim we first note that part of the proof of Theorem
\ref{prop} can be copied to show that the sets
$C_0,C_{\beta_\ell},\ldots,C_{(\kappa-1)\beta_\ell}$ are almost disjoint,
while, for $i\geq \kappa$, the sets $C_{i\beta_\ell}$ and
$C_{(i-\kappa)\beta_\ell}$ are equivalent with respect to $\mathcal{X}^{(m)}$.
Since $B_0$ is a Blackwell component for $\mathcal{X}^{(\ell)}$ and
$C_0\subset B_0$, we have $\cup_{i=0}^{k-1}C_{i\beta_\ell}\subset B_0$. But,
again in analogy with part of the proof of Theorem \ref{prop}, it is
easily seen that $\cup_{i=0}^{\kappa-1}C_{i\beta_\ell}$ is almost closed with
respect to $\mathcal{X}^{(\ell)}$, so, $B_0$ being atomic, we actually have
$\cup_{i=0}^{\kappa-1}C_{i\beta_\ell}\eqas B_0$. As in the proof of Theorem
\ref{prop} it is no restriction to assume that the sets
$C_0,C_{\beta_\ell},\ldots,C_{(\kappa-1)\beta_\ell}$ are disjoint rather than
almost disjoint.

Assuming that the Blackwell components for $\mathcal{X}^{(\ell)}$ are suitably
numbered, we have $C_1\subset B_1$ and the preceding argument can be
repeated to show that the sets
$C_1,C_{\beta_\ell +1},\ldots,C_{(\kappa-1)\beta_\ell +1}$ are disjoint, while
$\cup_{i=0}^{\kappa-1}C_{i\beta_\ell+1}\eqas B_1$.
Thus proceeding it follows eventually that
$\{C_0,C_1,\ldots,C_{\kappa\beta_\ell -1}\}$ constitutes a Blackwell
decomposition of $S$ for $\mathcal{X}^{(m)}$, so that
$\beta(\mathcal{X}^{(m)}) = \kappa\beta_\ell$, as claimed.

We finally observe that the $\kappa$ sets
$\cup_{i=0}^{\beta_\ell-1}C_{i\kappa+j},~ j=0,1,\dots,\kappa-1$,
are almost closed with respect to $\mathcal{X}^{(k)}$, so that $\kappa$
must be a divisor of $\beta(\mathcal{X}^{(k)})$. $\Box$ 
\bigskip

\noindent
This lemma has some interesting and useful corollaries, of which the first is
immediate.

\begin{corollary}
\label{cor2}
Let $\mathcal{X}$ be simple and $m> 1$. If $\ell$ is a divisor of $m$, then
$\beta(\mathcal{X}^{(\ell)})$ is a divisor of $\beta(\mathcal{X}^{(m)})$. 
\end{corollary}

\begin{corollary}
\label{cor3}
Let $\mathcal{X}$ be simple and $m>1$. If $\beta(\mathcal{X}^{(m)}) = m$, then
$\beta(\mathcal{X}^{(\ell)}) = \ell$ for all divisors $\ell$ of $m$.
\end{corollary}

\noindent
{\bf Proof.}~
Let $m=k\ell$. Then, by Lemma \ref{divisor},
\[
\beta(\mathcal{X}^{(m)})=k\ell=\kappa\beta(\mathcal{X}^{(\ell)}),
\]
with $\kappa$ a divisor of $\beta(\mathcal{X}^{(k)})$, and, hence, by Theorem
\ref{prop}, of $k$. Since $\beta(\mathcal{X}^{(\ell)})$ is a divisor of $\ell$ we
must have $\kappa=k$ and $\beta(\mathcal{X}^{(\ell)})=\ell$. $\Box$

\begin{corollary}
\label{cor4}
Let $\mathcal{X}$ be simple and $k,~\ell\geq 1$. Then
$\beta(\mathcal{X}^{(k\ell)}) =
\beta(\mathcal{X}^{(k)})\beta(\mathcal{X}^{(\ell)})$ if
$\beta(\mathcal{X}^{(k)})$ and $\beta(\mathcal{X}^{(\ell)})$ are relatively
prime. 
\end{corollary}

\noindent
{\bf Proof.}~
By Lemma \ref{divisor} we have 
$\beta(\mathcal{X}^{(k\ell)}) = \kappa \beta(\mathcal{X}^{(\ell)}) =
\lambda \beta(\mathcal{X}^{(k)})$, with $\kappa$ a divisor of
$\beta(\mathcal{X}^{(k)})$ and $\lambda$ a divisor of
$\beta(\mathcal{X}^{(\ell)})$. But if $\beta(\mathcal{X}^{(k)})$ and
$\beta(\mathcal{X}^{(\ell)})$ are relatively prime this is possible only if
$\kappa = \beta(\mathcal{X}^{(k)})$ and
$\lambda = \beta(\mathcal{X}^{(\ell)})$. $\Box$

%______________________________________________________________________________
\section{Asymptotic period}
\label{AP}

We are now ready to formally define the {\em asymptotic period\/} of a {\em
simple\/} Markov chain. As in the previous section, $\mathcal{X}$ denotes the
Markov chain of Section \ref{I}, and is, accordingly, stochastic, irreducible,
and aperiodic.
\bigskip

\noindent
{\bf Definition}~
Let the Markov chain $\mathcal{X}$ be simple. The {\em asymptotic period\/}
of $\mathcal{X}$ is given by
\begin{equation}
d(\mathcal{X}) := \sup\{m \geq 1\,|\,\beta(\mathcal{X}^{(m)})=m\};
\end{equation}
$\mathcal{X}$ is {\em asymptotically aperiodic\/} if $d(\mathcal{X})=1$,
otherwise $\mathcal{X}$ is {\em asymptotically periodic\/} with asymptotic
period $d(\mathcal{X}) > 1$.
\bigskip

\noindent
We shall see that it is possible for $\mathcal{X}$ to have $d(\mathcal{X})=\infty$.

If, for some $m$, we would have $\beta \equiv \beta(\mathcal{X}^{(m)}) >
d(\mathcal{X})$, then, by Lemma \ref{multiple}, $\beta(\mathcal{X}^{(\beta)})
= \beta > d(\mathcal{X})$, which is a contradiction. So we actually have the
following result, which formalizes the intuitive concept of asymptotic period put
forward in the introduction.

\begin{theorem}
The asymptotic period of a simple Markov chain $\mathcal{X}$ satisfies
\begin{equation}
d(\mathcal{X}) = \sup\{\beta(\mathcal{X}^{(m)})\,|\,m\geq 1\}.
\end{equation}
\end{theorem}

\noindent
From Theorem \ref{prop} we immediately conclude the following.

\begin{theorem}
If $\mathcal{X}$ is simple and nonatomic then $\mathcal{X}$ is asymptotically
aperiodic.
\end{theorem}

\noindent
The next theorem confirms the intimation in the introduction that an asymptotic
period larger than one requires the chain to be transient.

\begin{theorem}
\label{rec}
If $\mathcal{X}$ is simple and recurrent then $\mathcal{X}$ is asymptotically
aperiodic.
\end{theorem}

\noindent
{\bf Proof.}~
Suppose $d(\mathcal{X})>1$, so that $\beta(\mathcal{X}^{(m)}) = m$ for
some $m > 1$. Let $\{C_0,C_1,\ldots,C_{m-1}\}$ be a cyclic Blackwell
decomposition for $\mathcal{X}^{(m)}$, and choose $i\in C_0$. As a
consequence of Theorem \ref{prop} and the recurrence of state $i$ we must
have $P^{(\ell m+1)}(i,C_1)=1$ for all $\ell$. On the other hand, the
aperiodicity of $i$ implies $P^{(\ell m+1)}(i,i)> 0$ for $\ell$ sufficiently
large, contradicting the fact that $C_0$ and $C_1$ are disjoint. So
$\mathcal{X}$ must be asymptotically aperiodic if it is recurrent.
$\Box$
\bigskip

\noindent
An example of a chain with an asymptotic period greater than $1$ is obtained
by letting $\mathcal{X}$ be a transient birth-death process on the nonnegative
integers (as defined in the introduction) with self-transition probabilities
$r_i = 0$ except $r_0 = 1 - p_0 > 0$. Clearly, $\mathcal{X}$ is irreducible and
aperiodic, while Lemma \ref{bounded} implies that $\mathcal{X}$ is simple
(and atomic). But it is readily seen that $\beta(\mathcal{X}^{(2)})=2$, so
that $d(\mathcal{X}) > 1$.
(We will see in the next section that, actually, $d(\mathcal{X}) = 2$.)

It is possible for the asymptotic period of a Markov chain to be infinity.
Indeed, let us assume that the birth probabilities $p_i$ in a birth-death
process are such that $\prod_{i=0}^\infty p_i > 0.$
Then there is a probability $\prod_{i=j}^\infty p_i \geq \prod_{i=0}^\infty p_i$
that a visit to state $j$ is followed solely by jumps to the right. Hence, with
probability one, the process will make only a finite number of self-transitions
or jumps to the left. It follows that the sets
$C_i := \{i,n+i,2n+i,\ldots\},~i=0,1,\ldots,n-1,$ are (disjoint) atomic almost
closed sets with respect to $\mathcal{X}^{(n)}$, so that
$\beta(\mathcal{X}^{(n)})=n$ for all $n$ and, hence, $d(\mathcal{X})=\infty.$

Some further conditions for a simple Markov chain to be asymptotically
aperiodic are given next.

\begin{theorem}
Let $\mathcal{X}$ be a simple Markov chain. Then the following are
equivalent:\\
$(i)$	 $\mathcal{X}$ is asymptotically aperiodic;\\
$(ii)$	 $\mathcal{X}^{(m)}$ is simple for all $m > 1$;\\
$(iii)$  $\mathcal{X}^{(m)}$ is simple for all prime numbers $m$.
\end{theorem}

\noindent
{\bf Proof.}~
By Corollary \ref{cor1} the first statement implies the second. Evidently, the
second statement implies the third. To show that the third statement implies
the first, suppose $\beta(\mathcal{X}^{(m)})=1$ for all primes $m$. If
$d \equiv d(\mathcal{X})>1$, then $\beta(\mathcal{X}^{(d)}) = d$ and $d$ must
have a prime factor $p>1$. But then, by Corollary \ref{cor4},
$\beta(\mathcal{X}^{(p)}) =p$, which is impossible. $\Box$
\bigskip

\noindent
It may be desirable to have an upper bound on the asymptotic period of a Markov
chain. The next theorem, involving the condition
\begin{equation}
\label{cond(n)}
\begin{array}{l}
\mbox{there exists a constant $\delta > 0$ such that
$P^{(n)}(i,i) \geq \delta$ for all but}\\
\mbox{finitely many states $i \in S$},
\end{array}
\end{equation}
provides a criterion which may be used for this purpose.

\begin{theorem}
\label{upper}
If, for some $n$, the simple Markov chain $\mathcal{X}$ satisfies
condition \eqref{cond(n)}, then $d(\mathcal{X})$ is a divisor of $~n$.
\end{theorem}

\noindent
{\bf Proof.}~
In view of Theorem \ref{rec} we may assume that $\mathcal{X}$ is transient.
Suppose $\beta(\mathcal{X}^{(m)}) = m$ for some $m \geq 1,$ and let
$\{C_0,C_1,\ldots,C_{m-1}\}$ be a cyclic Blackwell decomposition for
$\mathcal{X}^{(m)}$. If \eqref{cond(n)} holds, then $P^{(n)}(i,C_0) \ge \delta$
for all but finitely many states $i\in C_0$. As a consequence 
$\mathbb{P}(U^{(m)}(C_0) \,|\, \theta^n U^{(m)}(C_0)) = 1$, 
and hence,
\[
\mathbb{P}(L^{(m)}(C_0)\,|\,\theta^n L^{(m)}(C_0)) = 1,
\]
since $C_0$ is almost closed with respect to $\mathcal{X}^{(m)}$. But by
Theorem \ref{prop} this is possible only if
$C_0 = C_{n\,(\mbox{\scriptsize{mod}}\,m)}$, that is, if $m$ is a divisor
of $n$. The result follows by definition of $d(\mathcal{X})$.
$\Box$
\bigskip

\noindent
We conclude this section with two corollaries of Theorem \ref{upper}, the
first one being evident.

\begin{corollary}
\label{cond(1)}
If the simple Markov chain $\mathcal{X}$ is such that $P(i,i)\geq \delta$ for some
$\delta>0$ and all but finitely many states $i\in S$, then $\mathcal{X}$ is
asymptotically aperiodic.
\end{corollary}

\begin{corollary}
\label{upper1}
If, for some $n$, the simple Markov chain $\mathcal{X}$ satisfies condition
(\ref{cond(n)}) while $\mathcal{X}^{(n)}$ is simple, then $\mathcal{X}$ is
asymptotically aperiodic.
\end{corollary}

\noindent
{\bf Proof.}~
If $\mathcal{X}$ satisfies (\ref{cond(n)}), then, by Theorem \ref{upper},
$d \equiv d(\mathcal{X})$ is a divisor of $n$, so that, by Corollary
\ref{cor2}, $\beta(\mathcal{X}^{(d)})$ is a divisor of
$\beta(\mathcal{X}^{(n)})$. Hence we must have
$d = \beta(\mathcal{X}^{(d)}) = 1$ if $\beta(\mathcal{X}^{(n)}) = 1$. $\Box$

%______________________________________________________________________________
\section{Birth-death processes}
\label{BD}

Throughout this section $S=\{0,1,\dots\}$ and $\mathcal{X}$ is a stochastic
and irreducible birth-death process on $S$ with at least one positive
self-transition probability, so that $\mathcal{X}$ is aperiodic.
Note that $\mathcal{X}$ is simple and atomic by Lemma \ref{bounded}.
As before, $P$ denotes the matrix of 1-step transition probabilities of
$\mathcal{X}$, and we use the notation \eqref{bdrates}. Letting
\[
\pi_0 := 1,~~ \pi_n := \frac{p_0p_1\ldots p_{n-1}}{q_1q_2\ldots q_n},
\quad n\geq 1, 
\]
and
\begin{equation}
\label{KL}
K_n := \sum_{j=0}^n \pi_j, \quad L_n := \sum_{j=0}^n \frac{1}{p_j\pi_j}.
\quad 0\leq n\leq\infty,
\end{equation}
we observe that $K_\infty+L_\infty=\infty$, and recall that
\begin{equation}
\label{classification}
\mathcal{X} ~\mbox{is}~
\left\{
\begin{array}{l@{}l}
\mbox{positive~recurrent} &~~\llr~~ K_\infty < \infty,~~ L_\infty = \infty\\
\mbox{null~recurrent} &~~\llr~~ K_\infty = \infty,~~ L_\infty = \infty\\
\mbox{transient} &~~\llr~~K_\infty = \infty,~~ L_\infty < \infty.
\end{array}
\right.
\end{equation}

\begin{theorem}
\label{bdap}
The asymptotic period $d(\mathcal{X})$ of the birth-death process $\mathcal{X}$
equals $1$, $2$, or $\infty$.
\end{theorem}

\noindent
{\bf Proof.}~
Suppose $2 < d \equiv d(\mathcal{X}) < \infty$, and let $\{C_0,C_1,\dots,C_{d-1}\}$
be a cyclic Blackwell decomposition of $S$ for $\mathcal{X}^{(d)}$.

By \eqref{structure0} we have, for $\ell=0,1\dots$ and $k$ sufficiently large,
$X(k+\ell) \in C_{\ell \,(\mbox{\scriptsize{mod}}\,d)}$ if $X(k) \in C_0$, and
in particular $X(k+1)\in C_1$.
Since $C_0$ and $C_1$ are disjoint, $X(k+1)=X(k)$ is impossible, but also
$X(k+1)=X(k)-1$ leads to a contradiction. Indeed, if $X(k)\in C_0$
and $X(k+1)=X(k)-1\in C_1$ then $X(k+2)\in C_2$ and hence
$X(k+2)=X(k)-2$, since the other options would contradict the fact that
$C_0$, $C_1$ and $C_2$ are disjoint. Thus continuing we eventually find that
$X(k+X(k)-1)=1\in C_{X(k)-1 \,(\mbox{\scriptsize{mod}}\,d)}$ and
$X(k+X(k))=0\in C_{X(k) \,(\mbox{\scriptsize{mod}}\,d)}$. But this would
imply $X(k+X(k)+1)=0$ or $X(k+X(k)+1)=1$, which is impossible since
$C_{X(k)-1 \,(\mbox{\scriptsize{mod}}\,d)}$,
$C_{X(k) \,(\mbox{\scriptsize{mod}}\,d)}$ and
$C_{X(k)+1 \,(\mbox{\scriptsize{mod}}\,d)}$ are disjoint.

So, assuming $k$ sufficiently large and $X(k)\in C_0$, we must have
$X(k+1)=X(k)+1\in C_1$. Repeating the argument leads to the conclusion
that for $k$ sufficiently large, $X(k)\in C_0$ implies
$X(k+\ell) = s+\ell\in C_{\ell \,(\mbox{\scriptsize{mod}}\,d)}$ for all
$\ell=0,1,\dots$.
We conclude that in the long run $\mathcal{X}$ will solely make jumps to the
right, that is, the number of self-transitions or jumps to the left will be
finite.  But then, as we have observed in Section \ref{AP},
$\beta(\mathcal{X}^{(n)})=n$ for all $n$, since the sets
$C'_i := \{i,n+i,2n+i,\ldots\},~i=0,1,\ldots,n-1,$ are (disjoint) atomic almost
closed sets with respect to $\mathcal{X}^{(n)}$. Hence $d(\mathcal{X})=\infty,$
contradicting our assumption $d(\mathcal{X})<\infty$.
$\Box$
\bigskip

\noindent
In what follows we derive necessary and sufficient conditions for
$d(\mathcal{X}) = \infty$ and for $d(\mathcal{X}) = 1$ (that is, for
asymptotic aperiodicity) in terms of the 1-step transition probabilities of the
process $\mathcal{X}$. By the above theorem we must have $d(\mathcal{X}) = 2$ in
the cases not covered by these criteria. The next theorem tells us when
$d(\mathcal{X}) = \infty$.

\begin{theorem}
\label{bdinf}
The birth-death process $\mathcal{X}$ has asymptotic period $d(\mathcal{X})
= \infty$ if and only if ~$\Pi_{i=0}^\infty p_i > 0$.
\end{theorem}

\noindent
{\bf Proof.}~
It has been shown already in Section \ref{AP} that $d(\mathcal{X})=\infty$ if
$\Pi_{i=0}^\infty p_i > 0$, so it remains to prove the converse. So suppose
$d(\mathcal{X})=\infty$ and let $d>2$ be such that $\beta(\mathcal{X}^{(d)})
=d$. The argument used in the proof of Theorem \ref{bdap} can be copied to
conclude that, with probability one, $\mathcal{X}$ will, in the long run, solely
make jumps to the right, but this obviously implies $\Pi_{i=0}^\infty p_i > 0$.
$\Box$
\bigskip

\noindent
A criterion for asymptotic aperiodicity in terms of the 1-step transition
probabilities follows after having established the validity of three lemmas.
The first is the following.

\begin{lemma}
\label{lbd1}
$\mathcal{X}$ is asymptotically aperiodic if and only if $\mathcal{X}^{(2)}$ is simple.
\end{lemma}

\noindent
{\bf Proof.}~
If $\mathcal{X}$ is asymptotically aperiodic, then, by definition,
$\beta(\mathcal{X}^{(2)})<2$, and hence $\beta(\mathcal{X}^{(2)})=1$,
that is, $\mathcal{X}^{(2)}$ is simple.
On the other hand, if $\mathcal{X}$ is {\em not\/} asymptotically aperiodic
then $d(\mathcal{X}) = 2$ or $d(\mathcal{X}) =\infty$, which both imply
$\beta(\mathcal{X}^{(2)})=2$, that is, $\mathcal{X}^{(2)}$ is {\em not\/}
simple. 
$\Box$
\bigskip

\noindent
Note that, by Lemma \ref{bounded}, $\mathcal{X}^{(2)}$ will be atomic if
it is simple.

In the next lemma a necessary and sufficient condition for $\mathcal{X}^{(2)}$
to be simple is given in terms of the polynomials $Q_n,~n \geq 0,$ that are
uniquely determined by the 1-step transition probabilities of $\mathcal{X}$ via
the recurrence relation
\begin{equation}
\label{recQ}
\begin{array}{l}
xQ_n(x)=q_nQ_{n-1}(x)+r_nQ_n(x)+p_nQ_{n+1}(x),\quad n > 1,\\
Q_0(x)=1,\quad  p_0Q_1(x)=x-r_0.
\end{array}
\end{equation}
The result is mentioned already in \cite[p.~275]{DS95}, but for completeness' sake
we give its proof.
\begin{lemma}
\label{lbd2}
$\mathcal{X}^{(2)}$ is simple if and only if ~$|Q_n(-1)| \to \infty$
as $n\to\infty$.
\end{lemma}

\noindent
{\bf Proof.}~
Writing $Q(x):=(Q_0(x),Q_1(x),\dots)^T$ (where superscript $T$ denotes
transposition), the recurrence relation \eqref{recQ} may be succinctly represented
by
\begin{equation}
P Q(x) = xQ(x).
\end{equation}
It follows that
\begin{equation}
P^2 Q(x) = x^2Q(x),
\end{equation}
so that the vectors $Q(1)$ and $Q(-1)$ are two distinct solutions of the system
of equations
\begin{equation}
\label{P2}
P^2y = y.
\end{equation}
Moreover, $P^2$ being a pentadiagonal matrix, any solution to \eqref{P2}
must be a linear combination of $Q(1)$ and $Q(-1)$. It follows that the
constant function is the only bounded harmonic function for
$P^2$ if and only if $Q_n(-1)$ is unbounded. Since $|Q_n(-1)|$ is increasing
(see Karlin and McGregor \cite[p.~76]{KM59} and Lemma \ref{lbd3} below),
Theorem \ref{blackwell2} leads to the required result.
$\Box$
\bigskip

\noindent
The third lemma constitutes an extension of Karlin and McGregor's result on the
sequence $\{Q_n(-1)\}_n$ referred to in the proof of the previous lemma.

\begin{lemma}
\label{lbd3}
The sequence $\{(-1)^nQ_n(-1)\}_n$ is increasing, and strictly increasing for $n$
sufficiently large. Moreover,
\[
\lim_{n\to\infty}(-1)^nQ_n(-1) = \infty ~~\llr~~
\sum_{j=0}^\infty \frac{1}{p_j\pi_j} \sum_{k=0}^j r_k\pi_k =\infty.
\]
\end{lemma}

\noindent
{\bf Proof.}~
Writing $\bQ_n(x):= (-1)^nQ_n(x)$ the recurrence relation \eqref{recQ} implies
\[
\begin{array}{l}
p_n\pi_n(\bQ_{n+1}(x)-\bQ_n(x)) = \\
\hspace{1cm} p_{n-1}\pi_{n-1}(\bQ_n(x)-\bQ_{n-1}(x))+(2r_n-1-x)\pi_n \bQ_n(x), \quad n \geq 1,\\
p_0\pi_0(\bQ_1(x)-\bQ_0(x)) = (2r_0-1-x)\pi_0 \bQ_0(x),
\end{array}
\]
so that
\[
p_n\pi_n(\bQ_{n+1}(x)-\bQ_n(x)) = \sum_{k=0}^n (2r_k-1-x)\pi_k \bQ_k(x), \quad n \geq 0,
\]
and hence
\begin{equation}
\label{qnxbar}
\bQ_{n+1}(x) = 1 + \sum_{j=0}^n \frac{1}{p_j\pi_j} \sum_{k=0}^j(2r_k-1-x) \pi_k \bQ_k(x),
\quad n \geq 0.
\end{equation}
It follows in particular (as observed already by Karlin and McGregor \cite[p.~76]{KM59})
that
\begin{equation}
\label{bQn1}
\bQ_{n+1}(-1) = 1 + 2\sum_{j=0}^n \frac{1}{p_j\pi_j} \sum_{k=0}^j r_k\pi_k \bQ_k(-1),
\quad n \geq 0,
\end{equation}
and hence 
\begin{equation}
\label{bQn2}
\bQ_{n+1}(-1) = \bQ_n(-1) + \frac{2}{p_n\pi_n} \sum_{k=0}^n  r_k\pi_k \bQ_k(-1),
\quad n \geq 0.
\end{equation}
Since $\bQ_0(-1)=1$ while $r_k>0$ for at least one state $k$ by the aperiodicity of
$\mathcal{X}$, the first statement follows. So we have $\bQ_n(-1)\geq 1$,
which, in view of \eqref{bQn1}, implies the necessity in the second statement.
To prove the sufficiency we let
\[
\be_j :=  \frac{2}{p_j\pi_j} \sum_{k=0}^j r_k\pi_k, \quad j \geq 0,
\]
and assume that $\sum_j \be_j$ converges. By \eqref{bQn2} we then have
\[
\bQ_{n+1}(-1)  \leq \bQ_n(-1)(1 + \beta_n), \quad n \geq 0,
\]
since $\bQ_n(-1)$ is increasing in $n$. It follows that
\[
\bQ_{n+1}(-1) \leq \prod_{j=0}^n (1 + \be_j), \quad n \geq 0.
\]
But, as is well known, $\prod_j (1 + \be_j)$ and $\sum_j \be_j$ converge together,
so we must have $\lim_{n\to\infty}\bQ_n(-1) < \infty$. 
$\Box$
\bigskip

\noindent
The Lemmas \ref{lbd1} -- \ref{lbd3} give us a necessary and sufficient condition
for $\mathcal{X}$ to be asymptotically aperiodic in terms of the $1$-step
transition probabilities.

\begin{theorem}
\label{bd1}
The birth-death process $\mathcal{X}$ is asymptotically aperiodic if and only if
\begin{equation}
\label{bd1f}
\sum_{j=0}^\infty \frac{1}{p_j\pi_j} \sum_{k=0}^j r_k\pi_k =\infty.
\end{equation}
\end{theorem}
\bigskip

\noindent
Considering \eqref{classification} and the fact that $r_k>0$ for at least one
state $k$ by the aperiodicity of $\mathcal{X}$, we see that $\mathcal{X}$ is
asymptotically aperiodic if $\mathcal{X}$ is recurrent, as we had observed
already in the more general setting of Theorem \ref{rec}. Another simple
sufficient condition for asymptotic aperiodicity is obtained by noting that
\[
\sum_{j=0}^n \frac{1}{p_j\pi_j} \sum_{k=0}^j r_k\pi_k 
\geq \sum_{j=0}^n \frac{r_j}{p_j},
\]
so that $\mathcal{X}$ is asymptotically aperiodic if $\sum_{j=0}^\infty
r_j/p_j = \infty$. Note that the latter condition is substantially weaker than the
condition given, in a more general setting, in Corollary \ref{cond(1)}.

%_______________________________________________________________________________

\end{document}